\documentclass[11pt]{article}
\usepackage[cp1251]{inputenc}
\usepackage[ukrainian, russian, english]{babel}

\usepackage{amsfonts, amssymb,mathrsfs}                                              %

\begin{document}
\selectlanguage{russian} \thispagestyle{empty}
 \pagestyle{myheadings}
\setcounter{page}{0}

{\bf Title:} Approximation of functions of several variables by linear methods in the space $S^p$

\vskip 5mm

{\bf Authors:} Viktor V. Savchuk and Andriy L. Shidlich$^{*}$

\vskip 5mm

{\bf Affilation:}  Department of Theory of Functions,

Institute of Mathematics of Ukrainian National Academy of Sciences,

3, Tereshchenkivska st., Kiev, Ukraine, 01601.

\vskip 5mm

{\bf email:} savchuk@imath.kiev.ua; andy709@list.ru, shidlich@imath.kiev.ua

\vskip 5mm

{\bf Abstract:} In the spaces $S^p$ of functions of several variables, $2\pi$-periodic in each variable, we  study the approximative properties of operators $A^\vartriangle_{\varrho,r}$ and $P^\vartriangle_{\varrho,s}$, which generate two summation methods of multiple Fourier series on triangular regions. In particular, in the terms of approximation estimates of these operators,  we give a constructive description of classes of functions, whose generalized derivatives belong to the classes $S^pH_\omega$.

\vskip 5mm


\vskip 5mm


\newpage

\begin{center}

{\bf  Approximation of functions of several variables  

by   linear methods in the space $S^p$}

\vskip 5mm

V{\small IKTOR} {V.}  {S}{\small AVCHUK} and A{\small NDRIY} {L.}  {S}{\small HIDLICH}

\end{center}

\vskip 5mm

\begin{minipage}[t]{14.5cm}

{\small {\it Abstract.} In the spaces $S^p$ of functions of several variables, $2\pi$-periodic in each variable, we  study the approximative properties of operators $A^\vartriangle_{\varrho,r}$ and $P^\vartriangle_{\varrho,s}$,  which generate two summation methods of multiple Fourier series on triangular regions. In particular, in the terms of approximation estimates of these operators,  we give a constructive description of classes of functions, whose generalized derivatives belong to the classes $S^pH_\omega$.}

       \end{minipage}

\vskip 5mm

{\bf 1. Introduction.} Let $d$ be an integer, let $\mathbb R^d,~\mathbb R^d_+,~\mathbb Z^d$ be the sets of all ordered sets 
${\bf k}: = (k_1, \ldots, k_d)$ with $d$ real numbers, $d$  real non-negative numbers and $d$  integers,  respectively. Set
$\mathbb T^d:=[0,2\pi]^d:= \underbrace{[0,2\pi]\times\ldots\times[0,2\pi]}_d$.

Further, let $L_p:=L_p(\mathbb T^d),~1\le p<\infty$ be the space of all functions $f$, 
defined on  $\mathbb R^d$, $2\pi$-periodic in each variable such that
\[
\|f\|_{_{\scriptstyle L_p}}:=\left(\int_{\mathbb
T^d}|f|^pd\sigma\right)^{1/p}<\infty,
\]
where $\sigma$ is the Lebesgue measure on  $\mathbb T^d$, and let  $C=C(\mathbb
T^d)$ be the space of all continuous functions on ${\mathbb R}^d$, $2\pi$-periodic in each variable, with the norm
$\|f\|_{_{\scriptstyle C}}:=\max\limits_{{\bf x}\in\mathbb
T^d}|f({\bf x})|.$

Let us set
\[
({\bf k,x}):=\sum_{j=1}^dk_jx_j
\]
and
\[
\widehat f({\bf k}):=\int_{\mathbb T^d}f({\bf x})e^{-i({\bf
k,x})}d\sigma({\bf x}),\quad{\bf k}\in\mathbb Z^d.
\]

The space $S^p,~1\le p<\infty,$ [1, Chap. XI] (see also [2]) 
is the space of all functions $f\in L_1$ such that
\[
\|f\|_{_{\scriptstyle S^p}}:=\left(\sum_{{\bf k}\in\mathbb
Z^d}|\widehat f({\bf k})|^p\right)^{1/p}<\infty.
\]
Functions $f\in L_1$ and $g\in L_1$ are equal in the space
$S^p$, if $\|f-g\|_{_{\scriptstyle S^p}}=0.$

Set
\[
H_{\nu}(f)({\bf x}):=\sum_{|{\bf k}|_1=\nu}\widehat f({\bf
k})e^{i({\bf k,x})},~|{\bf
k}|_1:=\sum_{j=1}^d|k_j|,\quad\nu\in\mathbb N.
\]
Then Fourier series of the function 
$f\in L_1$ can be represented as
\[
S[f]({\bf x}):=\sum_{{\bf k}\in {\mathbb Z}^d} \widehat f({\bf k})e^{i({\bf
x,k})}=\sum_{\nu=0}^\infty H_{\nu}(f)({\bf x}).
\]
Proceeding from the last equation, we consider the following linear operators
$S^\vartriangle_{n},$ $\sigma^\vartriangle_n$,
$P^\vartriangle_{\varrho,s}$ and $A^\vartriangle_{\varrho,r}$,
defined on $L_1$:
\[
S^\vartriangle_{n}(f)({\bf x})=\sum_{\nu=0}^{n}H_{\nu}(f)({\bf
x}),~n=0,1,\ldots,
\]
\[
\sigma^\vartriangle_n(f)({\bf
x})=\frac{1}{n+1}\sum_{\nu=0}^{n}S^\vartriangle_{\nu}(f)({\bf
x})=\sum_{\nu=0}^{n}\left(1-\frac{\nu}{n+1}\right)H_{\nu}(f)({\bf
x}),~n\in\mathbb N,
\]
\[
P^\vartriangle_{\varrho,s}(f)({\bf
x})=\sum_{\nu=1}^{\infty}\varrho^{\nu^s}H_{\nu}(f)({\bf x}),\quad
s>0,~\varrho\in[0,1),
\]
and
$$
A^\vartriangle_{\varrho,r}(f)({\bf x})=S^\vartriangle_{r-1}(f)({\bf
x})+\sum_{\nu=r}^{\infty}\lambda_{\nu,r}H_{\nu}(f)({\bf x}),
\eqno(1)
$$
where
\[
\lambda_{\nu,r}:=\lambda_{\nu,r}(\varrho):=\sum_{k=0}^{r-1}\bigg(\begin{array}{*{20}c}\nu\hfill\\
k\hfill\\
\end{array}\bigg)(1-\varrho)^{k}\varrho^{\nu-k}
=
\sum_{k=0}^{r-1}\frac{(1-\varrho)^{k}}{k!}~\frac{d^k}{d\varrho^k}\varrho^{\nu},\quad
r\in\mathbb N,~\varrho\in[0,1).
\]

The expressions $S^\vartriangle_n(f)({\bf x})$, $\sigma^\vartriangle_n(f)({\bf x})$ and $P^\vartriangle_{\varrho,s}(f)({\bf x})$ are called the triangular  partial sum of the Fourier series,  the triangular Fejer sum  and the generalised triangular Abel--Poisson sum  of the function  $f$, respectively. The expression  $A^\vartriangle_{\varrho,r}(f)({\bf x})$ is called the triangular  Abel--Poisson--Taylor sum of the function $f$.

Note that the definition of the operator  $A^\vartriangle_{\varrho,r}$ is correct, because 
\[
\sum_{k=0}^{r-1}\bigg(\begin{array}{*{20}c}\nu\hfill\\
k\hfill\\
\end{array}\bigg)(1-\varrho)^{k}\varrho^{\nu-k}\le
rq^{\nu}\nu^{r-1},~\mbox{where}~q=\max(1-\varrho,\varrho),
\]
and hence, for any function $f\in L_1$ and for any $0<\varrho<1$, the series on the right-hand
side of  (1) is majorized by the convergent series
\[
r\sum_{\nu=r}^{\infty}q^{\nu}\nu^{r-1}\sum_{|{\bf
k}|_1=\nu}|\widehat f({\bf k})|.
\]

Now we explain the motives  for the choice of a title of operators
$A^\vartriangle_{\varrho,r}(f)$.

Recall that the Poisson integral of a 
function $f \in L_1$ is the function $P(f)$, 
defined on $[0,1)^d\times {\mathbb R}^d$ by the equality
\[ P(f)({\bf \varrho,x})=\int_{\mathbb T^d}f({\bf x+t})P({\bf
\varrho,t})d\sigma({\bf t}),
\]
where
\[
P({\bf
\varrho,t}):=\prod_{j=1}^d\frac{1-\varrho_j^2}{1-2\varrho_j\cos
t_j+\varrho_j^2},~\varrho_j\in[0,1),
\]
is the multiple Poisson kernel and ${\bf x+t}:=(x_1+t_1,\ldots,x_d+t_d).$

In what follows,  the expression $P(f)(\varrho,{\ bf x})$  
means the Poisson integral, where $\varrho$ is a vector with the same 
coordinates, i.e.  $\varrho=(\varrho,\ldots,\varrho).$

According to the decomposition of the Poisson kernel
in powers of $\varrho $, for any function $f\in L_1$, its Poisson integral $P(f)(\varrho,{\bf x}),$
with $\varrho\in[0,1),$  can be written in the form
\[
P(f)(\varrho,{\bf
x})=\sum_{\nu=0}^{\infty}\varrho^{\nu}H_{\nu}(f)({\bf x}).
\]
The sum of the right-hand side of this equality coincides with the sum of
the Abel--Poisson series $\sum_{\nu=0}^{\infty}H_{\nu}(f)({\bf x})$,
 or, what is the same, with the sum of  $P^\vartriangle_{\varrho,1}(f)({\bf x}).$
 For ${\bf x}={\bf 0}: = (0, \ldots, 0) $, we denote by $F(\varrho)$ the sum of this series and consider it as a function of the variable $\varrho.$ It is clear that the function $ F $ is analytic on $[0,1).$ 
 Therefore, in the neighborhood of $\varrho \in [0,1)$ for the functions $F$, the following 
 Taylor's formula is satisfied:
\[
F(t)=\sum_{k=0}^\infty\frac{F^{(k)}(\varrho)}{k!}(t-\varrho)^{k}.
\]
By direct computation we see that the partial sum of this series of order
$r-1$ for $t=1$ coincides with the sum $A^\vartriangle_{\varrho,r}(f)({\bf 0}).$
In particular, for $r = 1$, we obtain $F(\varrho)=A^\vartriangle_{\varrho,1}(f)({\bf
0})=P^\vartriangle_{\varrho,1}(f)({\bf 0}).$

Consequently, on the one hand, the sum of $A^\vartriangle_{\varrho,r}(f)(\bf 0)$ can be interpreted as the Taylor sum of order $r-1$ of the function $F$, and on the other hand, for $r =1$, it can be interpreted as the Abel--Poisson sum.

The purpose of this paper is to investigate the operators $A^\vartriangle_{\varrho,r}$ and $P^\vartriangle_{\varrho,s}$ as the linear methods of approximation of functions in the spaces $S^p.$
In this case, our attention is drawn to the relationship of the approximative properties of the sums $A^\vartriangle_{\varrho,r}(f)$
and $P^\vartriangle_{\varrho,s}(f)$ with the differential properties of the function $f$, namely the properties of the derivatives which are determineed as follows.

Let $\psi=\{\psi({\bf k})\}_{{\bf k}\in\mathbb
Z^d}$ be a multiple numerical sequence whose members are not all zero and
\[
\mathscr Z(\psi):=\mathscr Z^d(\psi):=\left\{{\bf k}\in\mathbb Z^d :
\psi({\bf k})=0\right\}.
\]
In what follows, assume that the number of elements of the set $ \mathscr Z (\psi) $ is finite.

If for the function $ f \in L_1 $, there exists the function $g \in L_1 $ such that
$$
S[f]({\bf x})=\sum_{{\bf k}\in\mathscr Z(\psi)}\widehat{f}({\bf
k})e^{i({\bf k,x})}+\sum_{{\bf k}\in\mathbb Z^d}\psi({\bf
k})\widehat{g}({\bf k})e^{i({\bf k,x})},
\eqno(2)
$$
then we say that for the function $f$, there exists $\psi$-derivative $g$, 
for which we use the notation $g=f^{\psi}.$ In this case,  if $\mathscr Z(\psi)=\varnothing$, then the first sum in (2) is set equal to zero.

It is clear that  for any function from the space $S^p$, its $\psi$-derivative is the unique up to the sum $\sum_{{\bf k}\in\mathscr Z(\psi)}a_{{\bf k}}e^{i({\bf k, x})}$, where $a_{{\bf k}}$ are any numbers.

This definition of $\psi$-derivative is adapted to the needs of the research described in this paper and it is not fundamentally different from the established concept of $\psi$-derivative of A.I.~Stepanets [1, Ch. XI, \S 2, 2].

In the paper, for functions from $L_1$,  we consider $\psi$-derivatives of the following two forms:
\[
1)\quad \psi({\bf k})=\nu^{-r},~\mbox{if}~|{\bf
k}|_1=\nu,~\nu=0,1,\ldots,~r\ge 0;
\]
\[
2)\quad \psi({\bf k})=\left\{\matrix{0,~&\mbox{if}~|{\bf
k}|_1=0,1,\ldots,r-1,\cr\cr
\displaystyle\frac{(\nu-r)!}{\nu!},~&\mbox{if}~|{\bf
k}|_1=\nu,~\nu\ge r,~r\in\mathbb N.}\right.
\]
In the first case, for $\psi$-derivative of $f$, we use the notation $f^{(r)}$ and in the second case, we use the notation $f_{[r]}$. If $r=0$  then we set $f^{(0)}=f_{[0]}=f.$

In the terms of Poisson integrals, one can give the following interpretation of the derivative $f_{[r]}$: 

Assume that  $\varrho\in[0,1)$, then
$$
P(f_{[r]})(\varrho,{\bf
x})=\varrho^r\frac{\partial^r}{\partial\varrho^r}P(f)(\varrho,{\bf
x})
\eqno(3)
$$
and by virtue of the well-known theorem on radial limit values of the Poisson integral (see, eg, [3]), for almost all
 ${\bf x}\in\mathbb R^d$
$$
f_{[r]}({\bf x})=\lim_{\varrho\to
1-}\frac{\partial^r}{\partial\varrho^r}P(f)(\varrho,{\bf x}).
$$
Also note that  $f^{(1)}=f_{[1]}$.

In general case,  the operators $P^\vartriangle_{\varrho,s}$ were  perhaps first considered as the aggregates of approximation of functions of one variable in [4, 5]. The operators $A^\vartriangle_{\varrho,r}$ were first studied in [6], where in the terms of these operators, the author gives the structural characteristic of Hardy-Lipschitz classes $H^r_p\mathop{\rm Lip}\alpha$  of one variable functions,  holomorphic in the unit circle in the complex plane. 
In special cases, when $r=s=1$, the operators $A^\vartriangle_{\varrho,1}$ and $P^\vartriangle_{\varrho,1}$ coincide with each other and generate the Abel--Poisson method of summation of multiple Fourier series in the triangular areas. The problem of approximation of  $2\pi$-periodic functions by Abel--Poisson sums  has a long history, full of many results. Here we mention only the books [7--9], which contain fundamental results in this subject.

\bigskip
{\bf 2. The main results.} In the formulation of the main results, we use the following notation. Assume that $\mathbb Z^d_+:=\mathbb
R^d_+\cap\mathbb Z^d$, $\mathbb Z^d_-:=\left\{{\bf k}\in\mathbb Z^d
: k_j<0, j=1,\ldots,d\right\},$ $Y:=\mathbb Z^d_+\cup\mathbb Z^d_-,$
and
\[
L_{p,Y}:=L_{p,Y}(\mathbb T^d):=\left\{f\in L_p: \widehat f({\bf
k})=0~\forall~{\bf k}\in\mathbb Z^d\setminus Y\right\}.
\]

Let $\omega$ be a function defined on the interval $[0,1]$. For the space $X$, where $X$ is one of the spaces $L_p,~S^{p}$ or $C,$ set
\[
XH_\omega:=\left\{f : f\in
X,~\|f-f_h\|_{_{\scriptstyle X}}=O(\omega(|h|)),~|h|\to 0\right\},
\]
where $f_h({\bf x}):=f({\bf x}+h),~{\bf
x}+h:=(x_1+h,\ldots,x_d+h),~h\in\mathbb R^1.$

Further, we consider the functions $\omega(t)$, $0\le t\le 1$, satisfying the following conditions 1)-4):

1) $\omega(t)$ is continuous on $[0,1]$;

2) $\omega(t)\uparrow$;

3) $\omega(t)\not=0$ for any $t\in (0,1]$;

4) $\omega(t)\to 0$ as $t\to 0$;

\noindent as well-known condition $({\mathscr B})$: $\displaystyle{\sum\limits_{v=n+1}^\infty \frac 1v\,\omega\bigg(\frac 1v\bigg)=O\bigg[\omega\bigg(\frac 1n\bigg)\bigg]}$ (see, eg. [10]).

 \bigskip

{\bf P{\small ROPOSITION } 1.} {\it Assume that $1\le p<\infty$, $f\in L_1(\mathbb T^d),~d\in\mathbb N$ and $\omega$ is the function, satisfying  conditions  1)--4) and $({\mathscr B})$. The following statements are equivalent:

1) $\left\|S^\vartriangle_n\left(f_{[1]}\right)\right\|_{_{\scriptstyle
S^p}}=O(n\omega(\frac 1n)),\quad n\to\infty;$

2) $\left\|f-\sigma^\vartriangle_n\left(f\right)\right\|_{_{\scriptstyle
S^p}}=O(\omega(\frac 1n)),\quad n\to\infty$.

\noindent Furthermore, if one of the statements 1) or 2) is satisfied, then 

3) $f\in S^pH_\omega$.

\noindent If $f\in L_{1,Y}(\mathbb T^d)$, then the statements 1)--3) are equivalent.}

Let us give some comments to Proposition 1. 

First, let us note that the implication $2)\Rightarrow 3)$
is the statement of the type direct and inverse theorem for Fejer methods [9].

For a given number $\alpha\in(0,1]$ and for the space $X$, where $X$ is one of the spaces $L_p,~S^{p}$ or $C,$ set
\[
\mathop{\rm Lip}(\alpha, X):=\left\{f : f\in
X,~\|f-f_h\|_{_{\scriptstyle X}}=O(1)|h|^{\alpha},~|h|\to 0\right\},
\]
where $f_h({\bf x}):=f({\bf x}+h),~{\bf
x}+h:=(x_1+h,\ldots,x_d+h),~h\in\mathbb R^1.$

In the papers [11-18], F.M\'{o}ricz investigates conditions of absolute convergence of Fourier series of functions of one and several variables. In particular, in [11], the author obtains conditions of absolute convergence of Fourier series of functions of one variable and conditions of belonging such functions to the classes $\mathop{\rm Lip}(\alpha, C)$, which are given  in the terms of estimates of rate of increase partial sums Fourier of their derivates. In [15], similar results were obtained for the functions of several variables  in the terms of rectangular partial Fourier sums of their mixed derivatives.  Proposition 1 is closely related to these results in the following way:

In the case, where $d=1$, $p=1$ and $\omega(t)=t^\alpha$, the implication $1)\Rightarrow 3)$ coincides with the statement $(i)$ of Theorem 1 [11]. And in the case, where $d>1$, the difference of the implication $1)\Rightarrow 3)$  from the similar results of the paper  [15]  is that we consider the triangular partial Fourier sums.

{\bf Proof.} 

$1)\Rightarrow 2).$  Set $a_{\nu}=\|H_{\nu}(f)\|_{_{\scriptstyle S^p}}$, $\nu=0,1,\ldots$. Then
$$
\left\|f-\sigma^\vartriangle_n(f)\right\|^p_{_{\scriptstyle
S^p}}=\frac{1}{(n+1)^p}\sum_{\nu=1}^n(\nu
a_{\nu})^p+\sum_{\nu=n+1}^{\infty}\frac{1}{\nu^p}(\nu a_{\nu})^p.
\eqno(4)
$$
For a fixed integer $N>n$,  applying the Abel transformation to the last sum in the expression

\[
\frac{1}{(n+1)^p}\sum_{\nu=1}^n(\nu
a_{\nu})^p+\sum_{\nu=n+1}^{N}\frac{1}{\nu^p}(\nu a_{\nu})^p,
\]
 we obtain 
$$
\frac{1}{(n+1)^p}\sum_{\nu=1}^n(\nu
a_{\nu})^p+\sum_{\nu=n+1}^{N}\frac{1}{\nu^p}(\nu a_{\nu})^p=
\frac{1}{(n+1)^p}\sum_{\nu=1}^n(\nu
a_{\nu})^p-\frac{1}{n^p}\sum_{\nu=1}^n(\nu
a_\nu)^p+\frac 1{N^p}\sum\limits_{k=1}^N (ka_k)^p+
$$
$$
+\sum_{\nu=n+1}^{N}\left(\frac{1}{(\nu-1)^p}-\frac{1}{\nu^p}\right)
\sum_{k=1}^{\nu-1}(ka_k)^p\le
p\sum_{\nu=n}^{N}\frac{1}{\nu^{p+1}}\left\|S^\vartriangle_\nu(f_{[1]})\right\|^p_{_{\scriptstyle
S^p}}+\frac 1{N^p}\sum\limits_{k=1}^N (ka_k)^p.
\eqno(5)
$$
By virtue of statement 1), the last sum in this relation tends to zero as $N\to \infty$. Since the function $\omega$ satisfies condition $({\mathscr B})$, then for any $N$
$$
 \sum_{\nu=n}^{N}\frac{1}{\nu^{p+1}}\left\|S^\vartriangle_\nu(f_{[1]})\right\|^p_{_{\scriptstyle
S^p}}\le O(1)\sum_{\nu=n}^{\infty}\frac{\omega^p(1/\nu)}{\nu}=O(\omega^p(1/\nu)),\quad n\to\infty.
\eqno(6)
$$
Combining relations (4)--(6), we conclude that indeed the statement  2) is true.

$2)\Rightarrow 1).$ We have
\[
\left\|S^\vartriangle_n\left(f_{[1]}\right)\right\|_{_{\scriptstyle
S^p}}=\sum_{\nu=1}^n(\nu
a_{\nu})^p\le (n+1)^p\left\|f-\sigma^\vartriangle_n(f)\right\|^p_{_{\scriptstyle
S^p}}=O(n\omega(1/n)),~n\to\infty.
\]

Therefore, if one of the statements  1) or 2) is true, then the other statement is also true. Furthermore, since $\omega(t)\to 0$ as $t\to 0$, then from statement 2), it follows that $f\in S^p$. Putting $n=1/[h],~h>0$, we get
\[
\|f-f_{h}\|^p_{_{\scriptstyle S^p}}=\sum_{\nu=0}^{\infty}a^p_{\nu}|1-e^{i\nu
h}|^p=2^p\sum_{\nu=1}^{\infty}a^p_{\nu}\left|\sin\frac{\nu
h}{2}\right|^p\le
\sum_{\nu=1}^{n-1}a^p_{\nu}\left|2\sin\frac{\nu
}{2n}\right|^p+2^p\sum_{\nu=n}^{\infty}a^p_{\nu}\le
\]
\[
\le
\frac{2^p}{n^p}\sum_{\nu=1}^{n-1}(\nu
a_{\nu})^p+2^p\sum_{\nu=n}^{\infty}a_{\nu}^p=2^p\left\|f-\sigma^\vartriangle_{
n-1}(f)\right\|^p_{_{\scriptstyle S^p}}
=O(\omega(1/n))=O(\omega(h)),\quad n\to\infty,
\]
and hence, $f\in S^pH_\omega$.

To complete the proof of Proposition 1, it remains to verify that  if $f\in L_{1,Y}(\mathbb T^d)$, then the implication $3)\Rightarrow 1)$ is true.


Putting $h_n:=\pi/n,~n\in\mathbb N,$  by virtue of the inequality
$\nu h_n\le\pi\sin(\nu h_n/2),$ which is valid for all $\nu=1,2,\ldots,n,$ we obtain 
\[
\|S^\vartriangle_n(f_{[1]})\|^p_{_{\scriptstyle
S^p}}=\sum_{\nu=1}^{n}\nu^pa^p_{\nu}\le
\frac{1}{h^p_n}\sum_{\nu=1}^n(\nu h_n)^pa^p_{\nu}
\le\frac{\pi}{h^p_n}\sum_{\nu=1}^na^p_{\nu}\left|\sin\frac{\nu
h_n}{2}\right|^p\le
\]
\[
\le\frac{\pi}{h^p_n}\sum_{\nu=1}^{\infty}a^p_{\nu}\left|\sin\frac{\nu
h_n}{2}\right|^p=O\bigg(\frac{\omega(h_n)}{h_n}\bigg)=O(n\omega(1/n)),\quad n\to\infty.
\]
Proposition 1 is proved.
\bigskip

In the following theorem, which is the main result of the paper, we give   the direct and inverse theorem of the appoximation of functions by the operator $A^\vartriangle_{\varrho,r}(\cdot)$ in the space $S^p$ in the terms of majorants $\omega$.

{\bf T{\small HEOREM } 1.} {\it Assume that $1\le p<\infty$, $r\in\mathbb N$, $f\in L_1(\mathbb T^d),~d\in\mathbb N$ and $\omega$ is the function, satisfying  conditions  1)--4) and $({\mathscr B})$. The following statements are equivalent:

1) $\|f-A^\vartriangle_{\varrho,r}(f)\|_{_{\scriptstyle
S^p}}=O((1-\varrho)^{r-1}\omega(1-\varrho)),\quad\varrho\to 1-;$

2) $\left\|P(f_{[r]})({\varrho},\cdot)\right\|_{_{\scriptstyle
S^p}}=O(\frac {\omega(1-\varrho)}{1-\varrho}),\quad\varrho\to 1-;$

\noindent Furthermore, if one of the statements 1) or 2) is satisfied, then 

3) $f_{[r-1]}\in S^pH_\omega$.

\noindent If $f\in L_{1,Y}(\mathbb T^d)$, then the statements 1)--3) are equivalent.}

\bigskip

Let us note that the implication $2)\Rightarrow 3)$ is the statement of the type Hardy-Littlewood theorems [19].

Consider the approximative properties of the sums
$P^\vartriangle_{\varrho,s}(f)$ in the space $S^p$.

Let us show that 
$$
\|f-P^\vartriangle_{\varrho,s}(f)\|^p_{_{\scriptstyle S^p}}\sim
\|f^{(s-1)}-P^\vartriangle_{\varrho,1}(f^{(s-1)})\|^p_{_{\scriptstyle
S^p}},\quad\varrho\to 1-.
\eqno(7)
$$
For this, we set $a_{\nu}:=\|H_{\nu}(f)\|_{_{\scriptstyle S^p}},~\nu=0,1,\ldots.$ Then
\[
\|f-P^\vartriangle_{\varrho,s}(f)\|^p_{_{\scriptstyle
S^p}}=\sum_{\nu=1}^{\infty}\left(1-\varrho^{\nu^s}\right)^pa^p_{\nu}.
\]
Since
\[
\lim_{\varrho\to
1-}\frac{1-\varrho^{\nu^s}}{1-\varrho}=\nu^{s},\quad \nu\in\mathbb
N,~s\ge 1,
\]
then 
\[
\lim_{\varrho \to
1-}\frac{\|f-P^\vartriangle_{\varrho,s}(f)\|^p_{_{\scriptstyle
S^p}}}{\|f^{(s-1)}-P^\vartriangle_{\varrho,1}(f^{(s-1)})\|^p_{_{\scriptstyle
S^p}}}=\lim_{\varrho \to
1-}\frac{\displaystyle\sum_{\nu=1}^{\infty}a^p_{\nu}\frac{\displaystyle
(1-\varrho^{\nu^s})^p}{\displaystyle
(1-\varrho)^p}}{\displaystyle\sum_{\nu=1}^{\infty}a^p_{\nu}\nu^{p(s-1)}\frac{\displaystyle
(1-\varrho^{\nu})^p}{\displaystyle (1-\varrho)^p}}=
\frac{\|f^{(s)}\|^p_{_{\scriptstyle
S^p}}}{\|f^{(s)}\|^p_{_{\scriptstyle S^p}}}=1,
\]
Hence, relation (7) is true.

It is clear that
\[
P^\vartriangle_{\varrho,1}(f)({\bf
x})=A^\vartriangle_{\varrho,1}(f)({\bf x}).
\]
Therefore, applying Theorem 1 to the function $f=g^{(s-1)}$ with  $r=1$ and taking into account 
relation (7), we obtain the following result.

\bigskip
{\bf T{\small HEOREM } 2.} {\it Assume that $1\le p<\infty$, $r\in\mathbb N$, $f\in L_1(\mathbb T^d),~d\in\mathbb N$ and $\omega$ is the function, satisfying  conditions  1)--4) and $({\mathscr B})$. The following statements are equivalent:

1) $\|f-P^\vartriangle_{\varrho,s}(f)\|_{_{\scriptstyle
S^p}}=O(\omega(1-\varrho)),\quad\varrho\to 1-;$

2) $\left\|P(f^{(s)})({\varrho},\cdot)\right\|_{_{\scriptstyle
S^p}}=O(\frac {\omega(1-\varrho)}{1-\varrho}),\quad\varrho\to 1-;$

\noindent Furthermore, if one of the statements 1) or 2) is satisfied, then 

3) $f^{(s-1)}\in S^pH_\omega$.

\noindent If $f\in L_{1,Y}(\mathbb T^d)$, then the statements 1)--3) are equivalent.}

\bigskip

Let us note that in the case, where $\omega(t)=t^\alpha$, $\alpha>0$, Theorems 1 and 2 were proved in  [20].

\bigskip

R{\small EMARK.} 1. For $d=1$, the space $L_{1,Y}(\mathbb T^1)$ coincides with the space $L_{1}(\mathbb T^1)$, and therefore the statements
1)--3) in Theorems 1 and 2 are equivalent without any reservations.

\bigskip
{\bf 3. Proof of the results.} It is shown above that the Theorem 2 follows from Theorem 1. Therefore, it remains to prove the truth of Theorem 1.

$1)\Rightarrow 2).$ Let, as previously, $a_{\nu}:=\|H_{\nu}(f)\|_{_{\scriptstyle S^p}},~\nu=0,1,\ldots.$ Then
\[
\|f\|_{_{\scriptstyle
S^p}}=\left(\sum_{\nu=0}^{\infty}\|H_{\nu}(f)\|^p_{_{\scriptstyle
S^p}}\right)^{1/p}=\left(\sum_{\nu=0}^{\infty}a^p_{\nu}\right)^{1/p}.
\]

Since
$$
\sum_{k=0}^{\nu}\left(\begin{array}{*{20}c}\nu\hfill\\
k\hfill\\
\end{array}\right)(1-\varrho)^k\varrho^{\nu-k}=\big((1-\varrho)+\varrho\big)^{\nu}=1,~\nu=0,1,\ldots,
\eqno(8)
$$
then
\[
\|f-A^\vartriangle_{\varrho,r}(f)\|^p_{_{\scriptstyle
S^p}}=\sum_{\nu=r}^{\infty}|1-\lambda_{\nu,r}(\varrho)|^pa_{\nu}^p=
\sum_{\nu=r}^{\infty}\left(\sum_{k=r}^{\nu}\left(\begin{array}{*{20}c}\nu\hfill\\
k\hfill\\
\end{array}\right)(1-\varrho)^k\varrho^{\nu-k}\right)^pa_{\nu}^p\ge
\]
\[
\ge(1-\varrho)^{rp}\sum_{\nu=r}^{\infty}\left(\begin{array}{*{20}c}\nu\hfill\\
r\hfill\\
\end{array}\right)^p\varrho^{(\nu-r)p}a_{\nu}^p.
\]
On the other hand
\[
\frac{1}{(r!)^p}\left\|\frac{\partial^r}{\partial\varrho^r}P(f)(\varrho,\cdot)\right\|^p_{_{\scriptstyle
S^p}}=\sum_{\nu=r}^{\infty}\left(\begin{array}{*{20}c}\nu\hfill\\
r\hfill\\
\end{array}\right)^pa^p_{\nu}\varrho^{(\nu-r)p}.
\]
According to these relations and equality (3), we see that for $\varrho\to 1-$,
\[
\left\|P(f_{[r]})(\varrho,\cdot)\right\|_{_{\scriptstyle S^p}}\le
r!(1-\varrho)^{-r}\|f-A^\vartriangle_{\varrho,r}(f)\|_{_{\scriptstyle
S^p}}=O\bigg(\frac {\omega(1-\varrho)}{1-\varrho}\bigg).
\]

Further, for any numbers $n>r$ and $\varrho\in[0,1)$,
\[
\frac{1}{(r!)^p}\left\|\frac{\partial^r}{\partial\varrho^r}P(f)(\varrho,\cdot)\right\|^p_{_{\scriptstyle
S^p}}=\sum_{\nu=r}^{\infty}\left(\begin{array}{*{20}c}\nu\hfill\\
r\hfill\\
\end{array}\right)^pa^p_{\nu}\varrho^{(\nu-r)p}\ge
\]
\[
\ge
\varrho^{(n-r)p}\sum_{\nu=r}^{n}\left(\begin{array}{*{20}c}\nu\hfill\\
r\hfill\\
\end{array}\right)^pa^p_{\nu}=
\varrho^{(n-r)p}\frac{1}{(r!)^p}\left\|S^\vartriangle_n\left(f_{[r]}\right)\right\|^p_{_{\scriptstyle
S^p}}.
\]
From the last relation, putting  $\varrho=1-1/n$ and taking into account condition 2), we see that 
\[
\left\|S^\vartriangle_n\left(f_{[r]}\right)\right\|_{_{\scriptstyle
S^p}}\le O(1)\left(1-\frac{1}{n}\right)^{-n}
\frac {\omega(1/n)}{1/n}=O(n\omega(1/n)),\ \mbox{\rm as}\ n\to\infty. 
\]
Hence, if one of the statements 1) or 2) is true, then for the function  $g=f_{[r-1]}$, the statement 1) of Proposition 1  is also true:
$$
\left\|S^\vartriangle_n\left(g_{[1]}\right)\right\|_{_{\scriptstyle
S^p}}=O(n\omega(1/n)),\quad n\to\infty.
\eqno(9)
$$
According to Proposition 1, we also conclude that 
$$
\left\|g-\sigma^\vartriangle_n\left(g\right)\right\|_{_{\scriptstyle
S^p}}=\left\|f_{[r-1]}-\sigma^\vartriangle_n\left(f_{[r-1]}\right)\right\|_{_{\scriptstyle
S^p}}=O(\omega(1/n)),\quad n\to\infty,
\eqno(10)
$$
and the function $g=f_{[r-1]}$ as well as the function $f$ (by virtue of definition of derivative $f_{[r-1]}$) belongs to $S^pH_\omega$.
 
Let us verify the  validity of the implication $2)\Rightarrow 1).$

From identity (8), it follows that for any $\varrho\in[0,1]$,
\[
\sum_{k=r}^{\nu}\left(\begin{array}{*{20}c}\nu\hfill\\
k\hfill\\
\end{array}\right)(1-\varrho)^k\varrho^{\nu-k}\le 1,\quad\nu\ge r.
\]
This implies the relation
\[
\|f-A^\vartriangle_{\varrho,r}(f)\|^p_{_{\scriptstyle
S^p}}=\sum_{\nu=r}^{\infty}\left(\sum_{k=r}^{\nu}\left(\begin{array}{*{20}c}\nu\hfill\\
k\hfill\\
\end{array}\right)(1-\varrho)^k\varrho^{\nu-k}\right)^pa_{\nu}^p
\le\|f\|^p_{_{\scriptstyle S^p}}<\infty,
\]
From this relation, we conclude that for any
$\varepsilon>0$ there exists the number $n_0$ such that for all $n>n_0$ and all $\varrho\in[0,1)$,
$$
\|f-A^\vartriangle_{\varrho,r}(f)\|^p_{_{\scriptstyle
S^p}}=\sum_{\nu=r}^{n}\left(\sum_{k=r}^{\nu}\left(\begin{array}{*{20}c}\nu\hfill\\
k\hfill\\
\end{array}\right)(1-\varrho)^k\varrho^{\nu-k}\right)^pa_{\nu}^p+\varepsilon.
\eqno(11)
$$

Now, let us show that for all $\nu\ge r$, the following inequality is valid:
$$
\sum_{k=r}^{\nu}\left(\begin{array}{*{20}c}\nu\hfill\\
k\hfill\\
\end{array}\right)(1-\varrho)^k\varrho^{\nu-k}\le\left(\begin{array}{*{20}c}\nu\hfill\\
r\hfill\\
\end{array}\right)(1-\varrho)^r\quad\forall~\varrho\in[0,1].
\eqno(12)
$$

Putting $m=\nu-r$ and
\[
c_k=\frac{\displaystyle \left(\begin{array}{c}\nu\\
k+r\\
\end{array}\right)}{\displaystyle \displaystyle\left(\begin{array}{*{20}c}\nu\hfill\\
r\hfill\\
\end{array}\right)},\quad k=0,1,\ldots,m,
\]
we see that inequality (12) is true if and only if
\[
\sum_{k=0}^mc_k(1-\varrho)^k\varrho^{m-k}\le 1
\quad\forall~\varrho\in[0,1].
\]
To verify the validity of the last inequality it is sufficient to note that
\[
c_k=\frac{\nu!}{(k+r)!(\nu-k-r)!}\cdot\frac{r!(\nu-r)!}{\nu!}\le\frac{(\nu-r)!}{k!(\nu-r-k)!}=
\left(\begin{array}{c}m\\
k\\
\end{array}\right)
\]
and to use  binomial formula (see (8)).

Thus, by continuing further estimate (11), 
taking into account (12), we obtain
\[
\|f-A^\vartriangle_{\varrho,r}(f)\|^p_{_{\scriptstyle
S^p}}
\le
(1-\varrho)^{pr}\sum_{\nu=r}^{n}\left(\begin{array}{*{20}c}\nu\hfill\\
k\hfill\\
\end{array}\right)^pa_{\nu}^p+\varepsilon=
\frac{(1-\varrho)^{pr}}{(r!)^p}\|S^\vartriangle_{n}(f_{[r]})\|^p_{_{\scriptstyle
S^p}}+\varepsilon\le
\]
\[
\le
\frac{(1-\varrho)^{pr}(n+1)^p}{(r!)^p}\left\|f_{[r-1]}-\sigma^\vartriangle_n\left(f_{[r-1]}\right)\right\|_{_{\scriptstyle
S^p}}^p+\varepsilon.
\]
Now, let us set in these relations
$n=n_{\varrho}=[(1-\varrho)^{-1}]$, where $[\cdot]$ means the integer part of the number.  Then in view of (10), we get
\[
\|f-A^\vartriangle_{\varrho,r}(f)\|^p_{_{\scriptstyle
S^p}}=O(1)(1-\varrho)^{pr}n^p_{\varrho}\left\|f_{[r-1]}-\sigma^\vartriangle_{n_{\varrho}}\left(f_{[r-1]}\right)\right\|_{_{\scriptstyle
S^p}}^p+\varepsilon=
\]
\[
=O(1)(1-\varrho)^{p(r-1)}\omega(1-\varrho)+\varepsilon.
\]
as $\varrho\to 1-$. By virtue of arbitrary $\varepsilon$, from this relation it follows that the
implication $2) \Rightarrow 1)$ is true.

Therefore, we proved that statements 1) and 2) are equivalent. Furthermore, we proved that statements 1) and 2) are equivalent to relations (9) and (10).

The validity of the equivalence of statements 1)--3), in the case, where $f\in L_{1,Y}(\mathbb T^d)$,  follows from 
Proposition 1.

Theorem 1 is proved.

\bigskip

\footnotesize
 {\small


\begin{itemize}
  \item[{[1]}]
 A.\,I. Stepanets, Methods of approximation theory.\ VSP, Leiden (2005),  919~pp.
  \item[{[2]}]
A.I. Stepanets, Approximation characteristics of the spaces $S^p_\varphi$ in various metrics,\ Ukr. Mat. Zh., {\bf 53}, no. 8  (2001),
1121--1146; translation in Ukrainian Math. J., {\bf 53}, no. 8  (2001), 1340--1374.
  \item[{[3]}]
W. Rudin, Function theory in polydiscs.\ W. A. Benjamin Inc., New York-Amsterdam (1969),  188~pp.
  \item[{[4]}]
 Ja. S. Bugrov,  Bernstein type inequalities and their application to the study of differential properties of solutions of differential equations of higher order, Mathematica(Cluj), {\bf 5}, no.28 (1968), 5--25.
  \item[{[5]}]
 Ja. S. Bugrov,  Properties of solutions of differential equations of higher order in terms of weight classes,\ [in Russian] Studies in the theory of differentiable functions of several variables and its applications, IV. Trudy Mat. Inst. Steklov., {\bf 117} (1972), 47--61.
  \item[{[6]}]
V.V. Savchuk, Approximation of holomorphic functions by Taylor-Abel-Poisson means,\ Ukr. Mat. Zh., {\bf 59}, no. 9  (2007),
1253--1260; translation in Ukrainian Math. J., 59 (2007), no. 9, 1397--1407.
  \item[{[7]}]
N. K. Bari, Trigonometricheskie ryady  (Russian) [Trigonometric series]
 with the editorial collaboration of P. L. Ul'janov. Gosudarstv. Izdat. Fiz.-Mat. Lit., Moscow (1961) 936~pp.
  \item[{[8]}]
A. Zygmund, Trigonometricheskie ryady  (Russian) [Trigonometric series]: Vols. I, II,
Izdat. Mir, Moscow (1965) Vol. I: 615~pp.; Vol. II: 537 pp.
 \item[{[9]}]
P. Butzer and J. R. Nessel, Fourier analysis and approximation.\ Birh\"{a}user, Basel (1971),  553~pp.
\item[{[10]}]
N. K. Bari and S. B. Stechkin,  Best approximations and differential properties of two conjugate functions,\ Trudy Moskov. Mat. Obsch., 
{\bf 5} (1956), 483--522.
\item[{[11]}] 
F. M\'{o}ricz, Absolutely convergent Fourier series and function classes, J. Math. Anal. Appl., {\bf 324}, no. 2 (2006), 1168–1177. 
\item[{[12]}] 
F. M\'{o}ricz, Absolutely convergent Fourier series and generalized Lipschitz classes of functions, Colloq. Math., {\bf 113}, no. 1 (2008), 105–117.
\item[{[13]}] 
F. M\'{o}ricz, Absolutely convergent Fourier series and function classes. II, J. Math. Anal. Appl., {\bf 342}, no. 2 (2008), 1246–1249.
\item[{[14]}] 
F. M\'{o}ricz, Higher order Lipschitz classes of functions and absolutely convergent Fourier series, Acta Math. Hungar., {\bf 120}, no. 4 (2008), 355–366.
\item[{[15]}] 
F. M\'{o}ricz, Absolutely convergent multiple Fourier series and multiplicative Lipschitz classes of functions, Acta Math. Hungar., {\bf 121}, no. 1-2 (2008), 1–19.
\item[{[16]}] 
F. M\'{o}ricz, Absolutely convergent Fourier series, classical function classes and Paley's theorem, Anal. Math., {\bf 34}, no. 4 (2008), 261–276.
\item[{[17]}] 
F. M\'{o}ricz, Absolutely convergent Fourier series, enlarged Lipschitz and Zygmund classes of functions. East J. Approx., {\bf 15}, no. 1 (2009), 71–85.
\item[{[18]}] 
F. M\'{o}ricz,  Z. S\'{a}f\'{a}r, Absolutely convergent double Fourier series, enlarged Lipschitz and Zygmund classes of functions of two variables. East J. Approx., {\bf 16}, no. 1 (2010), 1–24.
\item[{[19]}] 
G.H. Hardy,  J.E. Littlewood,  Some properties of fractional integrals. II, Math. Z., {\bf 34}, no. 1 (1932), 403–439.
\item[{[20]}]  
V.V. Savchuk and A.L. Shidlich,  Approximation of several variables functions by linear methods in the space $S^p$,\ Zb. Pr. Inst. Mat. NAN Ukr., {\bf 4}, no. 1 (2007), 302--317.

\end{itemize}
\enddocument